\documentclass[a4paper,12pt]{article}

\usepackage{amsmath,amsfonts}
\usepackage{mathrsfs} 
\usepackage{parskip}
\usepackage{parskip}
\usepackage[T1]{fontenc}
\usepackage[utf8]{inputenc}
\usepackage{authblk}
\usepackage{enumerate}

\newtheorem{theorem}{Theorem}
\newtheorem{lemma}{Lemma}
\newtheorem{proposition}{Proposition}

\def\DD{\mathbb{D}}
\def\Pb{\mathbf{P}}

\def\Ex{\mathbf{E}}

\def\KK{\mathbb{K}}

\def\II{\mathbb{I}}

\def\UU{\mathbb{U}}

\def\sgn{{\rm sgn}} 
\def\1{\mbox{1\hspace{-.25em}I}}
\begin{document}
\title{On  Misspecifications in Regularity and Properties of Estimators}
\author[1]{O.V. Chernoyarov}
\author[2]{Yu.A. Kutoyants}
\author[3]{A.P. Trifonov}
\affil[1]{National Research University ``MPEI'', Moscow, Russia, }
\affil[2]{University of  Maine,  Le Mans,  France,}
\affil[2,3]{Voronezh State University, Voronezh, Russia}

\date{}
\maketitle
\

\begin{abstract}
 The problem of parameter estimation by the continuous time observations of a
 deterministic signal in white gaussian noise is considered.  The asymptotic
 properties of the maximul likelihood estimator are described in the
 asymptotics of small noise (large siglal-to-noise ratio). We are interested
 by the situation when there is a  misspecification in the regularity conditions. In
 particular it is supposed that the statistician uses a discontinuous
 (change-point type) model of signal, when the true signal is continuously
 differentiable function of the unknown parameter. 
\end{abstract}

\section{Introduction}

Consider the problem of parameter estimation by the observations of the
signals in White Gaussian Noise (WGN)
\begin{align*}
{\rm d}X_t=S\left(\vartheta ,t\right){\rm d}t+\varepsilon {\rm d}W_t,\quad
X_0=0,\quad 0\leq t\leq T.
\end{align*}
Here $S\left(\vartheta ,t\right)$ is deterministic known signal and we have
to estimate the parameter $\vartheta \in \Theta =\left(\alpha ,\beta \right)$
by continuous time observations $X^T=\left(X_t,0\leq t\leq T\right)$. We are
interested by the asymptotic behavior of estimators of this parameter in the
asymptotics of {\it small noise}, i.e.; as $\varepsilon \rightarrow 0$. It is
known that if the signal $S\left(\vartheta ,\cdot \right)$ is a smooth
function of $\vartheta $ with finite Fisher information 
\begin{align}
\label{0}
\II\left(\vartheta \right)=\int_{0}^{T}\dot S\left(\vartheta ,t\right)^2\,{\rm d}t,
\end{align}
then the maximim likelihood estimator $\hat\vartheta _\varepsilon $ is
consistent, asymptotically normal with the rate of convergence $\varepsilon $
and asymptotically efficient \cite{IH1}, \cite{T84}, \cite{TSh86}.  
Here in the sequel dot means derivation w.r.t. $\vartheta $. 

The situation changes if the signal $S\left(\vartheta
,t\right)=S\left(t-\vartheta \right)$, where $S\left(t\right)$  is a
discontinuous function of $t$, say, has a jump at the point $t=0$. Then the
Fisher information does not exist and the properties of the estimators are
essentially different. For example, the MLE  is consistent, has non gaussian limit 
distribution with the rate of convergence $\varepsilon ^2$ and asymptotically
efficient are bayesian estimators \cite{IH2}.

Let us recall that there is always a gap between mathematical model describing
the results of observations and the model which corresponds exactly to these
observations. Sometimes the difference is not important and the theoretical
results are in good agreement with the real data and sometimes the difference
can be essential. ``All models are wrong, but some are useful'' (G.E.P. Box).

We are interested by the situations, where the choosen models are not indeed
useful, i.e.; there is a misspecification. This misspecification concerns not
only the choice of the signal in the family of models close to the model of
real data, but we suppose that even the regularity conditions assumed by the
statistician are wrong. In particularly, the observed signal $S\left(\vartheta
,t\right)$ is {\it smooth}
with respect to the unknown parameter $\vartheta $, but the signal choosen by
the statistician $M\left(\vartheta,t\right)$ is discontinuous. Our goal is to
describe the properties of the corresponding pseudo-MLE. We remind  the
well-known property of this estimator that  it converges to the value
$\hat\vartheta $, which minimizes the Kullback-Leibler distance. Then we study
its limit distribution and show that it converges to non gaussian limit law
with the rate $\varepsilon ^{2/3}$. Remind that the real signals in
radiophysics can not have exactly rectangular form due to well-known physical
law and each time when the change-point type signal is used to describe the
real signal we are in situation of approximation \cite{T84}. This approximation can be
good or bad depending on the front of the signal and the level of
signal-to-noise ratio. 

We consider as well in some sense {\it inverse problem}, where the theoretical model is
smooth and the real data model is discontinuous and we describe the
asymptotics of the pseudo-MLE as $\varepsilon \rightarrow 0$. We show that in
this case the estimator $\hat\vartheta _\varepsilon $ converges to the point
$\hat\vartheta $ which minizes the Kulback-Leibler distance and is
asymptotically normal with the rate $\varepsilon $. 

At the end we describe the conditions on the misspecified model
(discontinuous vs discontinuous) which allow neveretheless to prove the 
consistency (true)  of the pseudo-MLE.

\section{Auxiliary results}

Let us consider the problem of parameter estimation by the observations (in
continuous time) of the deterministic signal in the presence of White Gaussian
Noise (WGN) of small intensity
\begin{align}
\label{01}
{\rm d}X_t=S\left(\vartheta ,t\right){\rm d}t+\varepsilon {\rm d}W_t,\quad
X_0=0,\quad 0\leq t\leq T, 
\end{align}
where the unknown parametr $\vartheta \in \Theta =\left(\alpha ,\beta
\right)$. We suppose for the simplicity of exposition that this parameter is
one-dimensional and that $\alpha $ and $\beta $ are finite. We are interested
by the  behavior of the estimators of this parameter  in the asymptotics of
{\it small noise}, i.e., as $\varepsilon \rightarrow 0$. It is well known that
in the case of the smooth (w.r.t. $\vartheta $ signal $S\left(\vartheta
,t\right)$ the maximum likelihood estimator (MLE) $\hat\vartheta _\varepsilon
$ and bayesian estimator (BE) $\tilde\vartheta _\varepsilon $
are consistent, asymptotically normal 
\begin{align*}
\varepsilon ^{-1}\left(\hat\vartheta _\varepsilon-\vartheta
\right)\Longrightarrow {\cal N}\left(0,\II\left(\vartheta
\right)^{-1}\right),\qquad \varepsilon ^{-1}\left(\tilde\vartheta
_\varepsilon-\vartheta 
\right)\Longrightarrow {\cal N}\left(0,\II\left(\vartheta \right)^{-1}\right),
\end{align*}
we have the convergence of all polynomial moments: for any $p>0$
\begin{align*}
\lim_{\varepsilon \rightarrow 0} \Ex_\vartheta
\left|\frac{\hat\vartheta _\varepsilon-\vartheta}{\varepsilon } \right|^p= \II\left(\vartheta
\right)^{-\frac{p}{2}}\Ex \left|\zeta \right|^p,\qquad \lim_{\varepsilon
  \rightarrow 0} \Ex_\vartheta 
\left|\frac{\tilde\vartheta _\varepsilon-\vartheta}{\varepsilon } \right|^p= \II\left(\vartheta
\right)^{-\frac{p}{2}}\Ex \left|\zeta \right|^p
\end{align*}
and the both estimators are asymptotically efficient. Here $\II\left(\vartheta
\right)$ is the Fisher information \eqref{0}.

Suppose that  the signal $S\left(\vartheta ,t\right)$ has {\it cusp}-type singularity,
say, $S\left(\vartheta ,t\right)=\left|t-\theta \right|^\kappa $, where
$0<\alpha<\vartheta <\beta <T $ and $\kappa \in (0,\frac{1}{2})$. Then the
Fisher information is $\infty $ and we have  a singular problem of parameter
estimation. Introduce the Hurst parameter $H=\kappa +\frac{1}{2}$ and
double-side fractional Brownian motion (fBm) $W^H\left(u\right), u\in R$.  It can
be shown that the MLE and BE  are consistyent, have limit distributions
\begin{align*}
\varepsilon ^{-\frac{2}{H}}\left(\hat\vartheta _\varepsilon-\vartheta
\right)\Longrightarrow \hat\xi ,\qquad \varepsilon ^{-\frac{2}{H}}\left(\tilde\vartheta
_\varepsilon-\vartheta 
\right)\Longrightarrow \tilde\xi ,
\end{align*}
 the polynomial moments converge : for any $p>0$
\begin{align*}
\lim_{\varepsilon \rightarrow 0} \Ex_\vartheta
\left|\frac{\hat\vartheta _\varepsilon-\vartheta}{\varepsilon^{\frac{2}{H}} } \right|^p= \Ex \left|\hat\xi \right|^p,\qquad \lim_{\varepsilon
  \rightarrow 0} \Ex_\vartheta 
\left|\frac{\tilde\vartheta _\varepsilon-\vartheta}{\varepsilon ^{\frac{2}{H}}
} \right|^p= \Ex \left|\tilde \xi \right|^p
\end{align*}
and the BE are asymptotically efficient. Here the random variables $\hat\xi$
and $\tilde\xi$ are defined by the relations 
\begin{align*}
Z\left(\hat\xi\right) =\sup_{u\in R}Z\left(u\right),\qquad \tilde \xi
=\frac{\int_{}^{}uZ\left(u\right){\rm d}u}{\int_{}^{}Z\left(u\right){\rm d}u}
\end{align*}
and the process 
\begin{align*}
Z\left(u\right)=\exp \left\{\gamma _\vartheta W^H\left(u\right)-\frac{\gamma
  _\vartheta  ^2}{2}\left|u\right|^ {2H}\right\},\qquad u\in R,
\end{align*}
where $\gamma _\vartheta $ is some constant.  The proofs can be carried out
using the general methods developped in \cite{IH81} as it was done in 
 \cite{DK03} in similar problem. 

Suppose now that the function $S\left(\vartheta ,t\right)$ has discontinuity,
say, 
\begin{align*}
S\left(\vartheta ,t\right)=h\left(t\right)\1_{\left\{t<\vartheta
  \right\}}+g\left(t\right)\1_{\left\{t\geq \vartheta \right\}},\qquad 0<\alpha <\vartheta <\beta <T,
\end{align*}
where $h\left(t\right)\not=g\left(t\right)$ for $t\in \left(\alpha ,\beta
\right)$. 

 It can
be shown that the MLE and BE  are consistent, have limit distributions
\begin{align*}
\varepsilon ^{-2}\left(\hat\vartheta _\varepsilon-\vartheta
\right)\Longrightarrow \hat\eta ,\qquad \varepsilon ^{-2}\left(\tilde\vartheta
_\varepsilon-\vartheta 
\right)\Longrightarrow \tilde\eta ,
\end{align*}
 the polynomial moments converge : for any $p>0$
\begin{align*}
\lim_{\varepsilon \rightarrow 0} \Ex_\vartheta
\left|\frac{\hat\vartheta _\varepsilon-\vartheta}{\varepsilon^{2 }}
\right|^p= \Ex \left|\hat\eta \right|^p,\qquad \lim_{\varepsilon   \rightarrow 0} \Ex_\vartheta 
\left|\frac{\tilde\vartheta _\varepsilon-\vartheta}{\varepsilon ^{2}
} \right|^p= \Ex \left|\tilde \eta \right|^p
\end{align*}
and the BE are asymptotically efficient. Here the random variables $\hat\eta $
and $\tilde\eta $ are defined by the relations 
\begin{align*}
Z\left(\hat\eta \right) =\sup_{u\in R}Z\left(u\right),\qquad \tilde \eta 
=\frac{\int_{}^{}uZ\left(u\right){\rm d}u}{\int_{}^{}Z\left(u\right){\rm d}u}
\end{align*}
and the process 
\begin{align*}
Z\left(u\right)=\exp \left\{\delta \left(\vartheta\right) W\left(u\right)-\frac{\delta \left(\vartheta\right)  ^2}{2}\left|u\right|\right\},\qquad u\in R,
\end{align*}
where $\delta \left(\vartheta\right) =h\left(\vartheta \right)-g\left(\vartheta \right)$.

We are intrerested by the following problem of misspecification. Suppose that
model of observations the choosen by the statistician is 
\begin{align*}
{\rm d}X_t=M\left(\vartheta ,t\right){\rm d}t+\varepsilon {\rm d}W_t,\quad X_0=0,\quad 0\leq t\leq T
\end{align*}
 and the real data model \eqref{01} are different. Especially we study the
 situations where the regularity conditions of these models do not
 coincide. For example, the signal $S\left(\vartheta ,t\right)$ is a smooth function of
 $\vartheta $ (regular case), but the statistician supposes that the observed
 model has singularities of cusp or discontinuous types. 

We discuss the conditions of the consistency, the rates of convergence and the
limit distributions of the corresponding pseudo-MLE  in such situations.

\section{Main results}

We suppose that the observed  process ({\it real model})  is 
\begin{align}
\label{11}
{\rm d}X_t=S\left(\vartheta_0 ,t\right){\rm d}t+\varepsilon {\rm d}W_t,\quad
X_0=0,\quad 0\leq t\leq T, 
\end{align}
where $\vartheta _0$ is the true value  of unknown parameter and $S\left(\vartheta ,\cdot
\right)\in L_2\left(0,T\right)$. If we use the
{\it theoretical model}
\begin{align*}
{\rm d}X_t=M\left(\vartheta ,t\right){\rm d}t+\varepsilon {\rm d}W_t,\quad
X_0=0,\quad 0\leq t\leq T, 
\end{align*}
with $M\left(\vartheta ,\cdot \right)\in L_2\left(0,T\right)$, then the
likelihood ratio (misspecified) is 
\begin{align}
\label{12}
V\left(\vartheta ,X^T\right)=\exp\left\{\frac{1}{\varepsilon ^2}\int_{0}^{T}
M\left(\vartheta ,t \right) {\rm d}X_t-\frac{1}{2\varepsilon ^2}\int_{0}^{T}
M\left(\vartheta ,t \right)^2 {\rm d}t \right\},\; \vartheta \in \Theta 
\end{align}
and the (pseudo) MLE $\hat\vartheta _\varepsilon $ is defined by the equation
\begin{align}
\label{13}
V\left(\hat\vartheta _\varepsilon ,X^T\right)=\sup_{\vartheta \in \Theta
}V\left(\vartheta ,X^T\right). 
\end{align} 
To understand what is the limit of the MLE we write the likelihood ratio as
follows
\begin{align*}
&\varepsilon ^2\ln V\left(\vartheta ,X^T\right)\\
&\qquad=\varepsilon \int_{0}^{T}
M\left(\vartheta ,t \right) {\rm d}W_t -\frac{1}{2}\int_{0}^{T}
\left[M\left(\vartheta ,t \right)^2 -2M\left(\vartheta ,t
  \right)S\left(\vartheta _0,t\right)\right]{\rm d}t \\
&\qquad=\varepsilon \int_{0}^{T}
M\left(\vartheta ,t \right) {\rm d}W_t -\frac{1}{2}
\left\|M\left(\vartheta ,\cdot \right) -S\left(\vartheta _0,\cdot\right)
\right\|^2+\left\|S\left(\vartheta _0,\cdot \right) \right\|^2 
\end{align*}
where we denoted as $ \left\|\cdot \right\|$ the $L_2\left(0,T\right)$ norm. 
It can be easily verified that under mild regularity conditions we have the convergence
\begin{align*}
\sup_{\vartheta \in \Theta }\left|\varepsilon ^2\ln L\left(\vartheta ,X^T\right)- \frac{1}{2}
\left\|M\left(\vartheta ,\cdot \right) -S\left(\vartheta _0,\cdot\right)
\right\|^2+\left\|S\left(\vartheta _0,\cdot \right) \right\|^2  \right|\rightarrow 0.
\end{align*}
Hence if we suppose that the equation 
\begin{align*}
\inf_{\vartheta} \left\|M(\vartheta ,\cdot ) -S\left(\vartheta _0,\cdot\right)
\right\|=\left\|M(\hat\vartheta ,\cdot ) -S\left(\vartheta _0,\cdot\right)
\right\|
\end{align*}
has a unique solution $\hat\vartheta $, then we obtain the well-known result  that
in the case of misspecification the MLE $\hat\vartheta _\varepsilon $
converges to the value $\hat\vartheta $, which minimizes the Kullback-Leibler
distance. It is interesting to note that in general case $\hat\vartheta
\not=\vartheta _0$ but sometimes $\hat\vartheta =\vartheta
_0$ and we consider the conditions of the consistency in such situations. The
most interesting for us is the question of the rate of convergence of the MLE
to the true value.

\subsection{Discontinuous versus smooth}

Here we consider the situation where the true model (described by
$S\left(\vartheta ,\cdot \right)$ is smooth w.r.t. $\vartheta $ but the
theoretical model choosen by statistician is discontinuous. We start with a
simple example. This example allows to see how we can have different rates of
convergence of estimators and what is the limit distribution of these
estimators. 

{\bf Example 1.} Suppose that the observed process is
\begin{align}
\label{rm}
{\rm d}X_t=\left(t-\vartheta _0\right){\rm d}t+\varepsilon {\rm d}W_t,\quad
X_0=0,\quad 0\leq t\leq T ,
\end{align}
where $\vartheta_0 \in \Theta =\left(\alpha ,\beta \right)$, $0<\alpha <\beta <T$ and
the theoretical model is 
\begin{align}
\label{tm}
{\rm d}X_t=\sgn\left(t-\vartheta \right){\rm d}t+\varepsilon {\rm d}W_t,\quad
X_0=0,\quad 0\leq t\leq T ,\quad \vartheta \in \Theta .
\end{align}
The pseudo-likelihood ratio is the function
\begin{align*}
V\left(\vartheta ,X^T \right)=\exp\left\{\frac{1}{\varepsilon^2
}\int_{0}^{T} \sgn\left(t-\vartheta \right){\rm d}X_t-\frac{T}{2\varepsilon
  ^2}\right\},\qquad \vartheta \in \Theta
\end{align*}
because $\sgn\left(t-\vartheta \right)^2=1$. Note that
\begin{align*}
\hat\vartheta =\arg\inf_{\vartheta \in \Theta }\int_{0}^{T}
\left[\sgn\left(t-\vartheta \right)-\left(t-\vartheta_0 \right)\right]^2{\rm d}t =\vartheta _0.
\end{align*} 
Hence the MLE $\hat \vartheta _\varepsilon $ defined by the relation 
\begin{align*}
\hat \vartheta _\varepsilon =\arg\sup_{\vartheta \in\Theta  }\int_{0}^{T}
\sgn\left(t-\vartheta \right){\rm d}X_t 
\end{align*}
in this misspecified parameter estimation problem is consistent.

To study its rate of convergence and the limit distribution we introduce a
normalized likelihood ratio 
\begin{align*}
Z_\varepsilon \left(u\right)=\frac{V\left(\vartheta_u ,X^T\right)}{V\left(\vartheta _0 ,X^T \right)},\qquad u\in \UU_\varepsilon ,
\end{align*}
where $\vartheta_u=\vartheta _0 +\varphi _\varepsilon u$. Here $\varphi
_\varepsilon \rightarrow 0$ will be choosen later and 
\begin{align*}
 \UU_\varepsilon =\left(\frac{\alpha -\vartheta _0 }{\varphi _\varepsilon },\frac{\beta
  -\vartheta _0 }{\varphi _\varepsilon }\right) \longrightarrow \left(-\infty ,\infty \right)
\end{align*}
as $\varepsilon \rightarrow 0$.

The substitution of the observation process in the likelihood ratio yields us
the following expression (we suppose that $u>0$)
\begin{align}
\label{5}
\ln Z_\varepsilon \left(u\right)&=\frac{1}{\varepsilon
}\int_{0}^{T}\left[ \sgn\left(t-\vartheta_0-\varphi _\varepsilon u
  \right)-\sgn\left(t-\vartheta_0 \right)\right]{\rm d}W_t\nonumber\\
&\quad +\frac{1}{\varepsilon^2
}\int_{0}^{T}\left[ \sgn\left(t-\vartheta_0-\varphi _\varepsilon u
  \right)-\sgn\left(t-\vartheta_0 \right)\right]\left(t-\vartheta
_0\right){\rm d}t \nonumber\\
&=-\frac{2}{\varepsilon}\int_{0}^{T}\1_{\left\{\vartheta _0 <t< \vartheta _0
  +\varphi _\varepsilon u\right\}}{\rm d}W_t
-\frac{2}{\varepsilon^2}\int_{0}^{T}\1_{\left\{\vartheta _0 <t< \vartheta _0 
  +\varphi _\varepsilon u\right\}} \left(t-\vartheta _0\right)  {\rm d}t\nonumber\\
&=-\frac{2}{\varepsilon}\left[ W_{\vartheta _0 
  +\varphi _\varepsilon u}-W_{\vartheta _0  }\right]-\frac{2}{\varepsilon^2}\int_{\vartheta _0}^{\vartheta _0 
  +\varphi _\varepsilon u}\left(t-\vartheta _0\right){\rm d}t \nonumber\\
&=\frac{2\sqrt{\varphi _\varepsilon }}{\varepsilon}
W\left(u\right)-\frac{\varphi _\varepsilon ^2}{\varepsilon ^2}u^2
=\frac{2\sqrt{\varphi _\varepsilon }}{\varepsilon}
\left[W_+\left(u\right)- \frac{\varphi _\varepsilon ^{3/2}}{\varepsilon } \frac{u^2}{2}\right],
\end{align}
where we denoted the Wiener process 
$$
W_+\left(u\right)=\varphi _\varepsilon ^{-1/2} \left[ W_{\vartheta _0 
  +\varphi _\varepsilon u}-W_{\vartheta _0  }\right],  \qquad u\in
\left.\left[0,\frac{\beta -\vartheta _0}{\varepsilon }\right. \right).
$$
Therefore if we take $\varphi _\varepsilon =\varepsilon ^{2/3}$, then we can write
\begin{align*}
\hat Z_\varepsilon \left(u\right)=\left(Z_\varepsilon
\left(u\right)\right)^{\frac{\varepsilon ^{2/3}}{2}}=\exp\left\{
W_+\left(u\right)-\frac{u^2}{2}\right\},  \qquad u\in \left.\left[0,\frac{\beta -\vartheta _0}{\varepsilon }\right. \right). 
\end{align*}
For the negative $u$ we obtain a similar representation
\begin{align*}
\hat Z_\varepsilon \left(u\right)=\left(Z_\varepsilon
\left(u\right)\right)^{\frac{\varepsilon ^{2/3}}{2}}=\exp\left\{
W_-\left(-u\right)-\frac{u^2}{2}\right\},  \qquad u\in \left.\left(\frac{\alpha 
    -\vartheta _0}{\varepsilon },0\right]\right. ,
\end{align*}
where $W_-\left(u\right),u\geq 0$ is a Wiener process  independent of
$W_+\left(u\right),u\geq 0$. Hence if we denote $W\left(\cdot \right)$ a double-sided
Wiener process, then the pseudo-likelihood ratio is
\begin{align*}
\hat Z_\varepsilon \left(u\right)=\left(Z_\varepsilon
\left(u\right)\right)^{\frac{\varepsilon ^{2/3}}{2}}=\exp\left\{
W\left(u\right)-\frac{u^2}{2}\right\},  \qquad u\in \UU_\varepsilon  .
\end{align*}
Now the properties of the MLE $\hat\vartheta _\varepsilon $ follow from the
relations 
\begin{align}
\label{an}
&\Pb_{\vartheta _0}\left(\frac{\hat\vartheta _\varepsilon-\vartheta _0
}{\varepsilon ^{2/3}}<x\right)=\Pb_{\vartheta _0}\left(\hat\vartheta
  _\varepsilon< \vartheta _0+\varepsilon ^{2/3}x\right) \nonumber\\
&\qquad \qquad =\Pb_{\vartheta _0}\left\{ \sup_{\vartheta <\vartheta
    _0+\varepsilon ^{2/3}x} V\left(\vartheta ,X^T\right)>\sup_{\vartheta \geq \vartheta
    _0+\varepsilon ^{2/3}x} V\left(\vartheta ,X^T\right)   \right\} \nonumber\\
&\qquad \qquad =\Pb_{\vartheta _0}\left\{ \sup_{\vartheta <\vartheta
    _0+\varepsilon ^{2/3}x} \frac{V\left(\vartheta ,X^T\right)}{
    V\left(\vartheta_0 ,X^T \right)}>\sup_{\vartheta \geq \vartheta 
    _0+\varepsilon ^{2/3}x} \frac{V\left(\vartheta ,X^T \right)
  }{V\left(\vartheta_0 ,X^T \right)} \right\}\nonumber \\ 
&\qquad \qquad =\Pb_{\vartheta _0}\left\{ \sup_{u <x,u\in\UU_\varepsilon }
  Z_\varepsilon \left(u\right)>\sup_{u \geq x,u\in\UU_\varepsilon }
  Z_\varepsilon \left(u\right) \right\}\nonumber\\ 
&\qquad \qquad =\Pb_{\vartheta _0}\left\{ \sup_{u <x,u\in\UU_\varepsilon }
  \hat Z_\varepsilon \left(u\right)>\sup_{u \geq x,u\in\UU_\varepsilon }
 \hat Z_\varepsilon \left(u\right) \right\} =\Pb_{\vartheta _0}\left(\hat u_\varepsilon <x\right)  ,
\end{align}
where we denoted $\hat u_\varepsilon $ the solution of the following equation
\begin{align*}
W\left(\hat u_\varepsilon \right)-\frac{\hat u_\varepsilon
  ^2}{2}=\sup_{u\in\UU_\varepsilon }\left(W\left(u\right)-\frac{u^2}{2}\right).
\end{align*}
It can be shown (see below) that 
\begin{align*}
\Pb_{\vartheta _0}\left(\hat u_\varepsilon <x\right) \longrightarrow
\Pb_{\vartheta _0}\left(\hat u <x\right)  ,
\end{align*}
i.e.;
\begin{align*}
\hat u_\varepsilon \Longrightarrow \hat u=\arg\sup_{u\in R}\left(W\left(u\right)-\frac{u^2}{2}\right).
\end{align*}

Therefore we obtain the following
\begin{proposition}
\label{P1} The pseudo-MLE $\hat\vartheta _\varepsilon $ in this problem is consistent, converges in
distribution
\begin{align*}
\frac{\hat\vartheta _\varepsilon -\vartheta _0}{\varepsilon
  ^{2/3}}\Longrightarrow \hat u
\end{align*}
and the moments converge: for any $p>0$
\begin{align*}
\lim_{\varepsilon \rightarrow 0}\Ex_{\vartheta _0}\left|\frac{\hat\vartheta
  _\varepsilon -\vartheta _0}{\varepsilon   ^{2/3}}  \right|^p =\Ex \left|\hat
u\right|^p .
\end{align*}

\end{proposition}
The proof follows from more general result of the Theorem \ref{T1} below. 

{\bf Remark 1.} Choosing different smooth signals in the class
\begin{align*}
{\cal S}=\left\{S\left(t-\vartheta \right)=\sgn\left(t-\vartheta \right)\left|t-\vartheta
\right|^\kappa , \kappa >\frac{1}{2}   \right\}
\end{align*}
and the same {\it theoretical model} \eqref{tm} we can obtain different rates
of convergence of estimators. Indeed, let us fix some $\kappa \in (\frac{1}{2},\infty)
  $. Then the corresponding
calculations like \eqref{5} provides us the expression
\begin{align*}
\ln Z_\varepsilon \left(u\right)=\frac{2\sqrt{\varphi _\varepsilon
}}{\varepsilon }\left[ W_+\left(u\right)-\frac{\varphi _\varepsilon ^{\frac{1}{2}+\kappa
  }}{\varepsilon  }  \frac{u^{1+\kappa }}{\left(1+\kappa \right) }\right].
\end{align*}
Therefore if we put $ \varphi _\varepsilon =\varepsilon ^{\frac{2}{2\kappa
    +1}} $, then 
\begin{align*}
\hat Z_\varepsilon
\left(u\right)=\exp\left\{W\left(u\right)-\frac{\left|u\right|^{1+\kappa
}}{1+\kappa }\right\} ,\qquad u\in \UU_\varepsilon 
\end{align*}
and the pseudo-MLE  $\hat\vartheta _\varepsilon $ is consistent and satisfies the relations
\begin{align*}
\frac{\hat\vartheta _\varepsilon -\vartheta _0}{\varepsilon ^{\frac{2}{2\kappa
    +1}}}=\hat u_\varepsilon =\arg\sup_{u\in \UU_\varepsilon }
\left[W\left(u\right)-\frac{\left|u\right|^{1+\kappa 
}}{1+\kappa }\right]\Longrightarrow \hat u,
\end{align*} 
where 
\begin{align*}
\hat u=\arg\sup_{u\in R}
\left[W\left(u\right)-\frac{\left|u\right|^{1+\kappa 
}}{1+\kappa }\right]
\end{align*}

Therefore choosing different $\kappa >\frac{1}{2}$ we can obtain any rate
$\varepsilon ^\gamma, \gamma <1$ of 
convergence of pseudo-MLE :
\begin{align*}
\frac{\hat\vartheta _\varepsilon -\vartheta _0}{\varepsilon ^\gamma }\Longrightarrow \hat u.
\end{align*}

\bigskip

Return now to the general smooth  model of observations
\begin{align}
\label{1-6}
{\rm d}X_t=S\left(\vartheta _0,t\right){\rm d}t+\varepsilon {\rm d}W_t,\quad
X_0=0,\quad 0\leq t\leq T 
\end{align}
and the discontinuous theoretical model
\begin{equation}
\label{1-7}
{\rm d}X_t=M\left(\vartheta ,t\right){\rm d}t+\varepsilon {\rm d}W_t,\quad
X_0=0,\quad 0\leq t\leq T ,
\end{equation}
where the signal
\begin{align*}
M\left(\vartheta ,t\right)=h\left(t\right)\1_{\left\{t<\vartheta
  \right\}}+g\left(t\right)\1_{\left\{t\geq \vartheta \right\}} .
\end{align*}
The unknown parameter $\vartheta \in \Theta =\left(\alpha ,\beta \right)$ with
$0<\alpha <\beta <T$.  We observe a trajectory $X^\varepsilon
=\left(X_t,0\leq t\leq T\right)$ of the solution of the equation \eqref{1-6} and we want to
estimate $\vartheta _0$ supposing that the observed process is
\eqref{1-7}. Therefore we introduce the pseudo-likelihood ratio
\begin{align*}
V\left(\vartheta ,X^T\right)&=\exp\left\{\frac{1}{\varepsilon
  ^2}\int_{0}^{\vartheta }h\left(t\right){\rm d}X_t+\frac{1}{\varepsilon
  ^2}\int_{\vartheta}^{T }g\left(t\right){\rm d}X_t\right.\\
&\left.\qquad\qquad  - \frac{1}{2\varepsilon
  ^2}\int_{0}^{\vartheta }h\left(t\right)^2{\rm d}t- \frac{1}{2\varepsilon
  ^2}\int_{\vartheta}^{T }g\left(t\right)^2{\rm d}t \right\},\quad \vartheta
\in \Theta 
\end{align*}
and define the pseudo-MLE $\hat\vartheta _\varepsilon $ by the equation
\begin{align*}
V\left(\hat\vartheta _\varepsilon ,X^T\right)=\sup_{\vartheta \in
  \Theta }V\left(\vartheta ,X^T\right).
\end{align*}

Let us introduce the following notations:
\begin{align*}
\delta \left(t \right)&=h\left(t \right)-g\left(t\right),\\
 \Phi\left(\vartheta \right)&=\int_{0}^{T }\left[M\left(\vartheta
  ,t\right)-S\left(\vartheta _0,t\right)\right]^2{\rm d}t,\\
 \ddot \Phi
\left(\vartheta \right)&=2\left[h\left(\vartheta \right)-S\left(\vartheta
  _0,\vartheta \right)\right]\left[\dot h\left(\vartheta \right)-S'\left(\vartheta
  _0,\vartheta \right)\right]\\ &-2\left[g\left(\vartheta
  \right)-S\left(\vartheta _0,\vartheta \right)\right]\left[\dot g\left(\vartheta
  \right)-S'\left(\vartheta _0,\vartheta \right)\right] ,\\
 \gamma (\hat \vartheta )&=\frac{\ddot \Phi(\hat\vartheta ) }{2},\qquad \hat\vartheta \in
\Theta ,\\ 
\hat Z\left(u\right)&=\exp\left\{\delta (\hat\vartheta
)W\left(u\right)-\frac{\gamma (\hat\vartheta )}{2}u^2\right\} ,\qquad u\in
R,\\
Z^o\left(v\right)&=\exp\left\{w\left(v\right)-\frac{v^2}{2}\right\} ,\qquad
 v\in R\\
 \hat u&=\arg\sup_{u\in R} \left[\delta (\hat\vartheta
  )W\left(u\right)-\frac{\gamma (\hat\vartheta )}{2}u^2\right],\qquad \hat
 v=\arg\sup_{u\in R} \left[w\left(v\right)-\frac{v^2}{2}\right]. 
\end{align*}
Here dot  means differentiating w.r.t. $\vartheta $,  prime means
differentiating w.r.t. $t $,    $W\left(u\right),u\in 
R$ and $w\left(v\right),v\in 
R$ are double-sided Wiener processes. 

 Note that
\begin{align}
\label{b}
\hat u= \hat v \left(\frac{\delta(\hat\vartheta )}{\gamma
  (\hat\vartheta)}\right) ^{\frac{2}{3}}. 
\end{align}
Indeed, let us put $u=rv$. Then we can write
\begin{align*}
&\delta (\hat\vartheta   )W\left(u\right)-\frac{\gamma (\hat\vartheta
    )}{2}u^2=\sqrt{r}\delta (\hat\vartheta   )w\left(v\right)-\frac{\gamma
    (\hat\vartheta )r^2}{2}v^2\\
&\qquad \qquad =\sqrt{r}\delta (\hat\vartheta   )\left[w\left(v\right)-\frac{\gamma
    (\hat\vartheta )r^{\frac{3}{2}}}{2\delta (\hat\vartheta   )}v^2
    \right]=\sqrt{r}\delta (\hat\vartheta   )\left[w\left(v\right)-\frac{v^2
    }{2}\right] 
\end{align*}
if we put $r=\delta (\hat\vartheta   )^{\frac{2}{3}}\gamma (\hat\vartheta
)^{-\frac{2}{3}}  $. Here $w\left(v\right)=r^{-1/2} W\left(rv\right)$. This proves \eqref{b}.

{\it Conditions ${\cal M}$.}
\begin{enumerate}
\item {\it  $\inf_{t\in \Theta }\delta \left(t\right)>0$.

\item The equation
\begin{align*}
\int_{0}^{\hat\vartheta }\left[h\left(t\right)-S\left(\vartheta
  _0,t\right)\right]^2{\rm d}t&+\int_{\hat\vartheta }^{T
}\left[g\left(t\right)-S\left(\vartheta _0,t\right)\right]^2{\rm d}t
=\inf_{\vartheta \in \Theta }   \Phi \left(\vartheta \right)
\end{align*}
has a unique solution $\hat\vartheta =\hat\vartheta \left(\vartheta
_0\right)\in \Theta $. 

\item The functions $h\left(t\right),g\left(t\right)$ and
$S\left(\vartheta ,t\right)$ are continuously differentiable
w.r.t. $t\in \Theta $.
\item  $\inf_{\vartheta \in \Theta } \ddot \Phi\left(\vartheta \right)>0$,
}
\end{enumerate}
Note that as $\hat\vartheta $ is the point of minimum of the function $\Phi
\left(\vartheta \right)$ we have the equality
\begin{align}
\label{eq1}
\dot \Phi \left(\hat\vartheta \right)=\left[h(\hat\vartheta)-S(\vartheta
  _0,\hat\vartheta )\right]^2-\left[g(\hat\vartheta )-S(\vartheta
  _0,\hat\vartheta )\right]^2=0,
\end{align}
which is equivalent to 
\begin{align*}
\left(h(\hat\vartheta )-g(\hat\vartheta
)\right)\left[h(\hat\vartheta )+g(\hat\vartheta )
  -2S(\vartheta _0,\hat\vartheta )\right]=0.
\end{align*}
Hence the point $\hat\vartheta $ satisfies to the equality
\begin{align*}
S(\vartheta _0,\hat\vartheta )=\frac{ h(\hat\vartheta )+g(\hat\vartheta ) }{2}.
\end{align*}
Of course, this is a necessary condition only. The equation
\begin{align}
\label{mm}
S(\vartheta _0,t )=\frac{ h(t )+g(t ) }{2},\qquad \alpha <t<\beta 
\end{align}
can have many solutions corresponding to the local maximums and minimums of
the function $\Phi \left(t\right),t\in \Theta $. If the equation \eqref{mm}
has no solution, say, 
\begin{align}
\label{mm1}
S(\vartheta _0,t )<\frac{ h(t )+g(t ) }{2},\qquad \alpha <t<\beta, 
\end{align}
then $\hat \vartheta =\alpha $. Otherwise $\hat \vartheta =\beta $. In these
two cases the behavior of the estimator $\hat\vartheta _\varepsilon $ can be
studied as it was done in \cite{Kut94}, Section 2.8. If we have the equality
 \begin{align*}
S(\vartheta _0,t )=\frac{ h(t )+g(t ) }{2},\qquad \alpha <t<\beta, 
\end{align*}
then any point of the interval $\left(\alpha ,\beta \right)$ can be taken as
$\hat \vartheta $. We do not study here the properties of $\hat\vartheta
_\varepsilon $ in this situation and in the situation when the function $\Phi
\left(\vartheta \right) ,\alpha <\vartheta <\beta $ has two or more points of
minimum. Note that such study can be done by the same way as in \cite{Kut94},
Section 2.7. 

 The properties of the pseudo-MLE $\hat\vartheta _\varepsilon $ are described
 in the following theorem.

\begin{theorem}
\label{T1} Let the conditions ${\cal M}$ be fulfilled then the estimator
$\hat\vartheta _\varepsilon $  converges to the value $\hat\vartheta $,  
has the limit distribution  
\begin{align}
\label{asno}
\frac{\hat\vartheta _\varepsilon-\hat\vartheta }{\varepsilon
  ^{2/3}}\Longrightarrow \hat u,
\end{align}
and  for any $p>0$
\begin{align}
\label{conv}
\lim_{\varepsilon \rightarrow 0}\Ex_{\vartheta _0} \left|\frac{\hat\vartheta
  _\varepsilon-\hat\vartheta }{\varepsilon ^{2/3}}\right|^p = \Ex_{\vartheta
  _0}\left|\hat u\right|^p.
\end{align}

\end{theorem}
{\bf Proof.} As before we study the normalized pseudo-likelihood ratio process
\begin{align*}
Z_\varepsilon \left(u\right)=\frac{V\left(\hat\vartheta +\varphi _\varepsilon
  u,X^T\right)}{V\left(\hat\vartheta,X^T \right)} ,\qquad
u\in \UU_\varepsilon , 
\end{align*}
where $\varphi _\varepsilon =\varepsilon ^{2/3}$. We have ($u>0$)
\begin{align*}
\ln Z_\varepsilon \left(u\right)&=\frac{1}{\varepsilon ^2}\int_{0}^{T} \left[
  M\left(\hat\vartheta +\varphi _\varepsilon u,t\right)-M\left(\hat\vartheta
  ,t\right)\right]{\rm d}X_t\\
&\qquad -\frac{1}{2\varepsilon ^2}\int_{0}^{T}  \left[
  M\left(\hat\vartheta +\varphi _\varepsilon u,t\right)^2-M\left(\hat\vartheta
  ,t\right)^2\right]{\rm d}t\\
&=\frac{1}{\varepsilon
}\int_{\hat\vartheta}^{\hat\vartheta+\varphi _\varepsilon
  u}\left[h\left(t\right)-g\left(t\right)\right] {\rm d}W_t \\
&\qquad -\frac{1}{2\varepsilon^2
}\int_{\hat\vartheta}^{\hat\vartheta+\varphi _\varepsilon
  u}\left(\left[h\left(t\right)-S\left(\vartheta _0,t\right)\right]^2-
\left[g\left(t\right)-S\left(\vartheta _0,t\right)\right]^2\right){\rm d}t\\ 
&=\frac{\delta (\hat\vartheta)\sqrt{\varphi
    _\varepsilon } }{\varepsilon } \left[\frac{W_{\hat\vartheta+\varphi _\varepsilon
  u}-W_{\hat\vartheta}}{\sqrt{\varphi _\varepsilon }} \right]-\frac{ \Phi
  \left(\hat\vartheta+\varphi _\varepsilon 
  u \right)-\Phi \left(\hat\vartheta \right)}{2\varepsilon ^2} +o\left(1\right)\\
&=\frac{\sqrt{\varphi _\varepsilon }}{\varepsilon }\delta (\hat\vartheta
)W_+\left(u\right) -\frac{\varphi _\varepsilon ^2u^2}{4\varepsilon
  ^2}\ddot \Phi(\hat\vartheta ) +o\left(1\right)\\
&=\varepsilon^{-2/3} \left[\delta (\hat\vartheta
)W_+\left(u\right)-  \frac{\ddot \Phi(\hat\vartheta )
  }{2}\,\frac{u^2}{2} \right] +o\left(1\right).
\end{align*}
Here we introduced the Wiener process
\begin{align*}
W_+\left(u\right)=\frac{W_{\hat\vartheta+\varphi _\varepsilon
  u}-W_{\hat\vartheta}}{\sqrt{\varphi _\varepsilon }},\quad u\in \left.\left[0,\frac{\beta -\hat\vartheta }{\varepsilon ^{2/3}}\right.\right),
\end{align*}
and used   in the expansion of $\Phi \left(\hat\vartheta +\varphi _\varepsilon
u\right)$ the equality  \eqref{eq1}.

For the negative values $u<0$ we obtain the similar representation
\begin{align*}
\ln Z_\varepsilon \left(u\right)={\varepsilon^{-2/3} }\left[\delta (\hat\vartheta
)W_-\left(-u\right)-  \frac{\ddot \Phi(\hat\vartheta )
  }{2}\,\frac{u^2}{2} \right] +o\left(1\right)
\end{align*}
with independent Wiener process $W_-\left(u\right),u\geq 0$.

Introduce  the random process 
\begin{align*}
\hat Z_\varepsilon \left(u\right)=\left(Z_\varepsilon
\left(u\right)\right)^{{\varepsilon ^{2/3}}{}}
 =
\exp\left\{\delta (\hat\vartheta )W\left(u\right)-\frac{\gamma (\hat\vartheta )}{2}
u^2+o\left(1\right) \right\} ,\quad u\in \UU_\varepsilon .
\end{align*}
We define $Z_\varepsilon \left(u\right)$ lineary decresing to zero on the interval 
$$
\left[\frac{\beta -\hat\vartheta}{\varepsilon^{2/3} }, \frac{\beta -\hat\vartheta}{\varepsilon^{2/3}  }+1\right]
$$
and increasing from zero to $\hat Z_\varepsilon \left(\frac{\alpha
  -\hat\vartheta}{\varepsilon } \right)$ on the interval 
$$
\left[\frac{\alpha  -\hat\vartheta}{\varepsilon^{2/3}  }-1, \frac{\alpha
    -\hat\vartheta}{\varepsilon^{2/3}  }\right] .
$$
Further we put $\hat Z_\varepsilon \left(u\right)=0 $ for 
\begin{align*}
u\not\in \left[\frac{\alpha  -\hat\vartheta}{\varepsilon^{2/3}  }-1,
  \frac{\beta -\hat\vartheta}{\varepsilon^{2/3}  }+1\right] 
\end{align*}
Now the process $\hat Z_\varepsilon \left(u\right)$ is defined for all $u\in
R$. Note that this process is continuous with probability 1. 

Let us denote by $Q_\varepsilon $ the measure induced by this process in the
space ${\cal C}_0\left(R\right)
$ of continuous functions decreasing to zero at infinity. The
corresponding measurable space we denote as $\left({\cal C}_0\left(R\right),
{\cal B}\right)$, where ${\cal B}$ is Borelian $\sigma $-algebra. By $Q$ we
denote the measure of the limit process $\hat Z\left(\cdot \right)$.

From this representation we obtain immediately the first lemma.
\begin{lemma}
\label{L0}
We have the
convergence of
finite-dimensional distributions of $\hat Z_\varepsilon \left(\cdot \right)$:
for any set $u_1,\ldots,u_k $ and any 
$k=1,2,\ldots$ 
\begin{equation}
\label{f-d}
\left(\hat Z_\varepsilon \left(u_1\right),\ldots,\hat Z_\varepsilon
\left(u_k\right)\right)\Longrightarrow  \left(\hat Z \left(u_1\right),\ldots,\hat Z
\left(u_k\right)\right).
\end{equation}
This convergence is uniform in $\vartheta $ on compacts $\KK\subset \Theta
$. 
\end{lemma}

We need the following elementary estimate
\begin{lemma}
\label{L1} There exists a constant $\kappa >0$ such that
\begin{equation}
\label{ng}
\Phi \left(\vartheta \right)-\Phi(\hat\vartheta )\geq \kappa \left(\vartheta -\hat\vartheta \right)^2.
\end{equation}
\end{lemma}
{\bf Proof.} As the point $\hat\vartheta $ is a unique minimum of  the function
$\Phi \left(\vartheta \right)$,  we can write for any $\nu >0$
$$
m\left(\nu \right)=\inf_{\left|\vartheta -\hat\vartheta \right|>\nu }\Phi
\left(\vartheta \right)-\Phi (\hat\vartheta )>0. 
$$ 
Hence for $\left|\vartheta -\hat\vartheta \right|>\nu $
\begin{align*}
\Phi \left(\vartheta \right)-\Phi (\hat\vartheta )\geq m\left(\nu
\right)\geq m\left(\nu\right)\frac{\left(\vartheta -\hat\vartheta
  \right)^2}{\left(\beta -\alpha \right)^2} .
\end{align*}
Further, for the values $\left|\vartheta -\hat\vartheta \right|\leq \nu $ by
Taylor expansion we have
\begin{align*}
\Phi \left(\vartheta \right)-\Phi (\hat\vartheta )=\frac{1}{2}\ddot \Phi(\hat\vartheta ) \left(\vartheta -\hat\vartheta \right)^2\left(1+o\left(1\right)\right).
\end{align*}
Therefore for sufficiently small $\nu $ we can write 
\begin{align*}
\Phi \left(\vartheta \right)-\Phi (\hat\vartheta )\geq \frac{1}{4}\ddot \Phi(\hat\vartheta )\left(\vartheta -\hat\vartheta \right)^2. 
\end{align*}
Taking 
\begin{align*}
\kappa =\min \left(  \frac{m\left(\nu\right)}{\left(\beta -\alpha \right)^2}
,\frac{1}{4}\ddot \Phi(\hat\vartheta ) 
\right)
\end{align*}
we obtain \eqref{ng}.

This estimate allows us to verify the boundness of all moments of the pseudo
likelihood ratio process.

\begin{lemma}
\label{L2}
For any $p>0$  there exist constants $c>0$ and $d>0$ such that for all
$\left|u\right|\geq d$ 
\begin{equation}
\label{tails}
\Ex_{\vartheta _0} \hat Z_\varepsilon ^{p}\left(u\right)\leq e^{-c u^2}.
\end{equation}
\end{lemma}
{\bf Proof.} Indeed, we have
\begin{align*}
\Ex_{\vartheta _0} \hat Z_\varepsilon
\left(u\right)^p=\exp\left\{\frac{p^2\varepsilon ^{-2/3}}{2}\int_{ 
  \hat\vartheta }^{\hat\vartheta+\varphi _\varepsilon u } \delta
\left(t\right)^2{\rm d}t-{p{\varepsilon ^{-4/3}}}{}
\left[\Phi(\hat\vartheta+\varphi _\varepsilon u) -\Phi(\hat\vartheta) \right]\right\}.
\end{align*}
Now the estimate \eqref{tails} follows from the relations 
\begin{align*}
&{\varepsilon ^{-2/3}}\int_{ \hat\vartheta }^{\hat\vartheta+\varphi _\varepsilon u
}\delta \left(t\right)^2{\rm d}t\leq {\sup_{t\in\Theta }\delta
  \left(t\right)^2}\left|u\right| ,\\ 
&{\varepsilon ^{-4/3}}{}
\left[\Phi(\hat\vartheta+\varphi _\varepsilon u) -\Phi(\hat\vartheta)
  \right]\geq {\kappa }u^2 ,
\end{align*}
where we used \eqref{ng}.  Therefore  we obtain the estimate \eqref{tails}
with some $c>0$ and $d>0$. 

\begin{lemma}
\label{L3} For $\left|u_1\right|<N,\left|u_2\right|<N$ and any $N>0$ we have
the estimate
\begin{equation}
\label{rr}
\Ex_{\vartheta _0}\left| \hat Z_\varepsilon \left({u_2}\right)-\hat
Z_\varepsilon \left({u_1}\right)\right|^4 \leq C\left(1+N^2\right)\left|u_2-u_1\right|^2
\end{equation}
with some constant $C>0$.
\end{lemma}
{\bf Proof.}
Let us denote
\begin{align*}
a_t&={\varepsilon ^{-1/3}\delta \left(t\right)},\qquad b_t=-{\varepsilon
  ^{-2/3} \delta \left(t\right)\left[
    h\left(t\right)+g\left(t\right)-2S\left(\vartheta _0,t\right)\right]}{},\\
G\left(t\right)&=\exp\left\{ \int_{\hat\vartheta +\varphi
   _\varepsilon u_1}^{t} a_s{\rm d}W_s+  \int_{\hat\vartheta +\varphi
   _\varepsilon u_1}^{t} b_s{\rm d}s \right\}
\end{align*}
Note that 
\begin{align*}
 G\left(\hat\vartheta +\varphi
   _\varepsilon u_2\right)=\frac{\hat Z_\varepsilon \left(u_2\right)}{\hat Z_\varepsilon
   \left(u_1\right)}.
\end{align*}
The process $G\left(t\right)$ has stochastic differential
\begin{align*}
{\rm d}G\left(t\right)=G\left(t\right)\left[{b_t}+\frac{a_t^2}{2}
  \right]{\rm d}t +G\left(t\right)a_t{\rm d}W_t,\quad G\left(\hat\vartheta +\varphi
   _\varepsilon u_1\right)=1.
\end{align*}
Therefore 
\begin{align*}
G\left(\hat\vartheta +\varphi _\varepsilon u_2\right)=1+\int_{\hat\vartheta
  +\varphi _\varepsilon u_1}^{\hat\vartheta +\varphi _\varepsilon
  u_2}G\left(t\right)\left[{b_t}+\frac{a_t^2}{2} \right]{\rm d}t+\int_{\hat\vartheta
  +\varphi _\varepsilon u_1}^{\hat\vartheta +\varphi _\varepsilon
  u_2}G\left(t\right)a_t{\rm d}W_t .
\end{align*}
We write
\begin{align*}
&\Ex_{\vartheta _0}\left| \hat Z_\varepsilon \left({u_2}\right)-\hat
Z_\varepsilon \left({u_1}\right)\right|^4 =\Ex_{\vartheta _0}\hat
Z_\varepsilon \left({u_1}\right)^4\left|G\left(\hat\vartheta +\varphi
_\varepsilon u_2\right)-1\right|^4\\
&\quad = \Ex_{\vartheta _0}\hat
Z_\varepsilon \left({u_1}\right)^4\left|  \int_{\hat\vartheta
  +\varphi _\varepsilon u_1}^{\hat\vartheta +\varphi _\varepsilon
  u_2}G\left(t\right)\left[{b_t}+\frac{a_t^2}{2} \right]{\rm d}t+\int_{\hat\vartheta
  +\varphi _\varepsilon u_1}^{\hat\vartheta +\varphi _\varepsilon
  u_2}G\left(t\right)a_t{\rm d}W_t      \right|^4\\
&\quad\leq  C_1\Ex_{\vartheta _0}\hat
Z_\varepsilon \left({u_1}\right)^4\left|  \int_{\hat\vartheta
  +\varphi _\varepsilon u_1}^{\hat\vartheta +\varphi _\varepsilon
  u_2}G\left(t\right)\left[{b_t}+\frac{a_t^2}{2} \right]{\rm d}t
\right|^4\\
&\qquad\qquad +C_2\Ex_{\vartheta _0}\hat 
Z_\varepsilon \left({u_1}\right)^4\left| \int_{\hat\vartheta 
  +\varphi _\varepsilon u_1}^{\hat\vartheta +\varphi _\varepsilon
  u_2}G\left(t\right)a_t{\rm d}W_t      \right|^4\\
&\quad\leq C_1 \left( u_2-u_1\right)^3\varphi _\varepsilon ^3  \int_{\hat\vartheta
  +\varphi _\varepsilon u_1}^{\hat\vartheta +\varphi _\varepsilon
  u_2}\Ex_{\vartheta _0}\hat
Z_\varepsilon \left({u_1}\right)^4G\left(t\right)^4
\left|{b_t}+\frac{a_t^2}{2}\right|^4  {\rm d}t \\
&\qquad\qquad +C_2\left(\Ex_{\vartheta _0}\hat 
Z_\varepsilon \left({u_1}\right)^8\Ex_{\vartheta _0}\left| \int_{\hat\vartheta 
  +\varphi _\varepsilon u_1}^{\hat\vartheta +\varphi _\varepsilon
  u_2}G\left(t\right)a_t{\rm d}W_t      \right|^8\right)^{1/2}.
\end{align*}
For stochastic integral we have the estimate
\begin{align*}
&\Ex_{\vartheta _0}\left| \int_{\hat\vartheta 
  +\varphi _\varepsilon u_1}^{\hat\vartheta +\varphi _\varepsilon
  u_2}G\left(t\right)a_t{\rm d}W_t      \right|^8\leq C\Ex_{\vartheta _0}\left( \int_{\hat\vartheta 
  +\varphi _\varepsilon u_1}^{\hat\vartheta +\varphi _\varepsilon
  u_2}G\left(t\right)^2a_t^2{\rm d}t\right)^4\\
&\qquad \qquad \leq \left( u_2-u_1\right)^3\varphi _\varepsilon ^3 \int_{\hat\vartheta 
  +\varphi _\varepsilon u_1}^{\hat\vartheta +\varphi _\varepsilon
  u_2}a_t^8\, \Ex_{\vartheta _0}G\left(t\right)^8\,{\rm d}t.
\end{align*}
Further
\begin{align*}
\Ex_{\vartheta _0}G\left(t\right)^8&=\exp\left\{ 32\int_{\hat\vartheta 
  +\varphi _\varepsilon u_1}^{t}a_s^2\; {\rm d}s-8\varepsilon
^{-4/3}\left[\Phi \left(t\right)-\Phi \left(\hat\vartheta  
  +\varphi _\varepsilon u_1\right)\right]\right\},\\
\Ex_{\vartheta _0}\hat Z_\varepsilon \left(u_1\right)^8&=\exp\left\{
32\int_{0}^{\hat\vartheta  
  +\varphi _\varepsilon u_1}a_s^2\; {\rm d}s-8\varepsilon ^{-4/3}\Phi \left(\hat\vartheta  
  +\varphi _\varepsilon u_1\right)\right\}
\end{align*}
Hence 
\begin{align*}
\Ex_{\vartheta _0}\left| \hat Z_\varepsilon \left({u_2}\right)-\hat
Z_\varepsilon \left({u_1}\right)\right|^4 &\leq
C\left|u_2-u_1\right|^4+C\left|u_2-u_1\right|^2\\
&\leq C\left(1+N^2\right)\left|u_2-u_1\right|^2 
\end{align*}
for $\left|u_1\right|\leq N$ and $\left|u_2\right|\leq N$.

Now the properties \eqref{asno} and \eqref{conv} of the pseudo MLE
$\hat\vartheta _\varepsilon $ follow from 
the Lemmae \ref{L0}, \ref{L2}, \ref{L3} and the Theorem 1.10.1 in \cite{IH81}. 

We see that the $\hat\vartheta _\varepsilon $ has a ``bad'' rate of
convergence. Note that for other estimators the rate can be better. 

Let us study the {\it trajectory fitting estimator} $\vartheta _\varepsilon
^*$ defined by the relation
\begin{align*}
\vartheta _\varepsilon^*=\arg\inf_{\vartheta \in \Theta } \int_{0}^{T}\left[
  X_t-m\left(\vartheta ,t\right)\right]^2{\rm d}t, 
\end{align*}
where 
\begin{align*}
m\left(\vartheta ,t\right)=\int_{0}^{t} M\left(\vartheta ,s\right)\,{\rm d}s.
\end{align*}
Suppose that the function
\begin{align*}
\Psi\left(\vartheta \right)=\int_{0}^{T}\left[m\left(\vartheta
  ,t\right)-s\left(\vartheta _0,t\right)\right]^2{\rm d}t,\qquad \vartheta \in \Theta  
\end{align*}
has a unique minimum at the point $\vartheta ^*\in \Theta $. Here 
\begin{align*}
s\left(\vartheta _0,t\right)=\int_{0}^{t}S\left(\vartheta _0,v\right)\,{\rm d}v.
\end{align*} 
This estimator admits the representation
\begin{align*}
\frac{\vartheta _\varepsilon ^*-\vartheta ^*}{\varepsilon }=\frac{ \int_{0}^{T}
W_t\,\dot m\left(\vartheta ^*,t\right){\rm d}t}{\int_{0}^{T}
\dot m\left(\vartheta ^*,t\right)^2{\rm d}t }\left(1+o\left(1\right)\right).
\end{align*}
Therefore this estimator is asymptotically normal with the rate $\varepsilon
$. The details of the proof can be found in the Section 7.4 in \cite{Kut94}, .

\subsection{Smooth versus discontinuous} 

Suppose now that the true model \eqref{1-6} has discontinuous trend
coefficient $S\left(\vartheta _0,t\right)$ of the following form 
\begin{align}
\label{20}
{\rm d}X_t=\left[h\left(t\right)\1_{\left\{t<\vartheta
  _0\right\}}+g\left(t\right)\1_{\left\{t\geq \vartheta _0\right\}}
  \right]{\rm d}t+\varepsilon {\rm d}W_t,\quad X_0=0,\quad 0\leq t\leq T ,
\end{align}
where  $\vartheta_0 \in \Theta =\left(\alpha ,\beta \right)$,
 $0<\alpha <\beta <T$, but the
statistician uses the  model
\begin{align}
\label{21}
{\rm d}X_t=M\left(\vartheta ,t\right){\rm d}t+\varepsilon {\rm d}W_t,\quad
X_0=0,\quad 0\leq t\leq T , 
\end{align}
with the ``smooth'' signal $M\left(\vartheta ,\cdot \right)$. The likelihood
ratio $L\left(\vartheta ,X^T\right)$ and the pseudo-MLE $\hat\vartheta
_\varepsilon $ are defined by the same relations \eqref{12}, 
\eqref{13}. As before,  we are interested by the asymptotic behavior of $\hat\vartheta
_\varepsilon $ as $\varepsilon \rightarrow 0$. 

To show that the situation is quite different we start with the example which
is ``symmetric'' to the Example 1. 

{\bf Example 2.} Suppose that the observed process is
\begin{align*}
{\rm d}X_t=\sgn\left(t-\vartheta _0\right){\rm d}t+\varepsilon {\rm
  d}W_t,\quad X_0=0,\quad 0\leq t\leq T , 
\end{align*}
and we use the model 
\begin{align*}
{\rm d}X_t=\left(t-\vartheta \right){\rm d}t+\varepsilon {\rm
  d}W_t,\quad X_0=0,\quad 0\leq t\leq T , 
\end{align*}
to estimate the parameter $\vartheta \in \Theta =\left(\alpha ,\beta \right)$,
where $0<\alpha <\beta <T$. 

It is easy to see that the function
\begin{align*}
\Phi \left(\vartheta \right)=\int_{0}^{T}\left[t-\vartheta
  -\sgn\left(t-\vartheta _0\right) \right]^2{\rm d}t ,\qquad \vartheta \in \Theta 
\end{align*}
atteints its minimum at the point
\begin{align*}
\hat\vartheta =\frac{T}{2}-1+\frac{2}{T}\;\vartheta _0
\end{align*}
and therefore the pseudo-MLE
\begin{align*}
\hat\vartheta _\varepsilon =\frac{T^2-2X_T}{2T}=\hat\vartheta
-\frac{W_T}{T}\,\varepsilon \longrightarrow  \hat\vartheta .
\end{align*}
It has Gaussian distribution 
\begin{align*}
\frac{\hat\vartheta _\varepsilon -\hat\vartheta }{\varepsilon }\sim {\cal
  N}\left(0 ,T^{-1}\right)
\end{align*} 
and the rate of convergence is $\varepsilon $. 

Note that if the true value is $\vartheta _0=\frac{T}{2}$, then it is
``consistent'' otherwise - not. Of course the consistent estimator can be
constructed as follows 
\begin{align*}
\vartheta _\varepsilon ^*=\left(\hat\vartheta
_\varepsilon-\frac{T}{2}+1\right)\frac{T}{2}\longrightarrow  \vartheta _0
\end{align*}
but to do this we need to know the true model. 

Let us return to the problem with the equations \eqref{20}, \eqref{21} and
introduce the conditions of regularity.

{\it Conditions ${\cal R}$.}
\begin{enumerate} 
\item {\it The functions $h\left(\cdot \right)$ and $g\left(\cdot \right)$ are
  bounded and for all $\vartheta _0\in \left[\alpha ,\beta \right]$ we have $
  h\left(\vartheta _0\right)\not=g\left(\vartheta _0\right)$. 
\item The function 
\begin{align*}
\Phi \left(\vartheta \right)=\int_{0}^{T}\left[M\left(\vartheta,t\right)
  -S\left(\vartheta _0,t\right) \right]^2{\rm d}t ,\qquad \vartheta \in \Theta 
\end{align*}
has a unique minimum at the point $\hat\vartheta\in\Theta  $.
\item The function $M\left(\vartheta ,t\right)$ is two times continuously
  differentiable w.r.t. $\vartheta $.
\item The function}
\begin{align*}
\ddot \Phi (\hat\vartheta )=2\int_{0}^{T}\ddot
M(\hat\vartheta,t)\left[M(\hat\vartheta,t) -S\left(\vartheta _0,t\right) \right]
  {\rm d}t +\int_{0}^{T}\dot M(\hat\vartheta,t)^2\,  {\rm d}t >0.
\end{align*}

\end{enumerate}

Let us denote
\begin{align*}
\II\left(\vartheta \right)=\int_{0}^{T}\dot 
M(\vartheta ,t)^2\,{\rm d}t ,\qquad   \DD\left(\vartheta _0\right)^2=\ddot \Phi(\hat\vartheta )^{-2}\II(\hat\vartheta ). 
\end{align*} 

\begin{theorem}
\label{T2} Let the conditions ${\cal R}$ be fulfilled, then the estimator
$\hat\vartheta _\varepsilon $ converges to the value $\hat\vartheta $, is
asymptotically normal
\begin{align*}
\frac{\hat\vartheta _\varepsilon -\hat\vartheta }{\varepsilon }\Longrightarrow \hat \xi \sim 
     {\cal N} \left(0, \DD\left(\vartheta _0\right)^2\right),
\end{align*}
and for any $p>0$
\begin{align*}
\lim_{\varepsilon \rightarrow 0}\Ex_{\vartheta _0}\left|\frac{\hat\vartheta
  _\varepsilon -\hat\vartheta }{\varepsilon }  \right|^p=\Ex_{\vartheta
  _0}\left|\hat \xi \right|^p. 
\end{align*}
\end{theorem}
{\bf Proof.} Using the Taylor expansion we can write for the likelihood ratio
\begin{align*}
Z_\varepsilon \left(u\right)=\frac{V\left(\hat\vartheta +\varepsilon
  u,X^T\right)}{V\left(\hat\vartheta ,X^T\right)} ,\quad u\in \UU_\varepsilon
=\left( \frac{\alpha -\hat\vartheta }{\varepsilon }, \frac{\beta
  -\hat\vartheta }{\varepsilon }\right) 
\end{align*}
the presentation
\begin{align*}
\ln Z_\varepsilon \left(u\right)=u\int_{0}^{T}\dot M\left(\hat\vartheta
,t\right){\rm d}W_t-\frac{u^2}{2} \ddot \Phi(\hat\vartheta )+ o\left(1\right).
\end{align*}
Therefore, if we denote 
\begin{align*}
Z\left(u\right)=\exp\left\{u\int_{0}^{T}\dot M\left(\hat\vartheta
,t\right){\rm d}W_t-\frac{u^2}{2} \ddot \Phi(\hat\vartheta ) \right\},\quad u\in R,
\end{align*}
then we obtain the first lemma.
\begin{lemma}
\label{L5}
We have the
convergence of
finite-dimensional distributions of $ Z_\varepsilon \left(\cdot \right)$:
for any set $u_1,\ldots,u_k $ and any 
$k=1,2,\ldots$ 
\begin{equation}
\label{f-d2}
\left( Z_\varepsilon \left(u_1\right),\ldots, Z_\varepsilon
\left(u_k\right)\right)\Longrightarrow  \left( Z \left(u_1\right),\ldots, Z
\left(u_k\right)\right).
\end{equation}
This convergence is uniforme in $\vartheta $ on compacts $\KK\subset \Theta $.
\end{lemma}
The next lemma can be proved following the same arguments as the Lemma
\ref{L1}. 

\begin{lemma}
\label{L6} There exists a constant $\kappa >0$ such that
\begin{equation}
\label{ng2}
\Phi \left(\vartheta \right)-\Phi(\hat\vartheta )\geq \kappa \left(\vartheta
-\hat\vartheta \right)^2.  
\end{equation}
\end{lemma}

Note that now the moments of $Z_\varepsilon \left(u\right)$ are no more
bounded. Denote $\vartheta _u=\hat\vartheta +\varepsilon u$, then for any
$\gamma >0$, we can write
\begin{align}
&\Ex_{\vartheta _0}Z_\varepsilon \left(u\right)^\gamma\nonumber\\
&\qquad  =\exp\left\{
\left({\frac{\gamma ^2}{2\varepsilon ^2}}\;\int_{0}^{T}\left[M(\vartheta
  _u)-M(\hat\vartheta )\right]^2{\rm d}t-\frac{\gamma }{2\varepsilon
  ^2}\;\left[\Phi(\vartheta_u )-\Phi(\hat\vartheta ) \right] \right)\right\}\nonumber\\
&\qquad \leq \exp\left\{
\left({\frac{\gamma ^2}{2}}\;M-\frac{\gamma\kappa  }{2}\right)\;  u^2
\right\}= 1,
\label{qq}
\end{align} 
where we denoted
\begin{align*}
M=\sup_{\vartheta \in \Theta }\int_{0}^{T}\dot M(\vartheta  )^2{\rm d}t
\end{align*}
and put $\gamma =M^{-1}\kappa $.

Therefore we introduce the following normalized likelihood ratio 
\begin{align*}
\hat Z_\varepsilon \left(u\right)=Z_\varepsilon \left(u\right)^\gamma ,\qquad
u\in \UU_\varepsilon 
\end{align*}
and establish the properties of this process similar to $\hat Z_\varepsilon
\left(\cdot \right)$ in  Lemmae \ref{L2} and \ref{L3}.
\begin{lemma}
\label{L7}
Suppose that the conditions ${\cal R}$ are fulfilled, then we have the
estimates 
\begin{align}
\label{qq1}
&\Ex_{\vartheta _0}\hat Z_\varepsilon ^{{1}/{2}}\left(u\right)\leq e^{-\kappa u^2},\\
&\Ex_{\vartheta _0}\left|\hat Z_\varepsilon^{{1}/{2}}
\left(u_2\right)-\hat Z_\varepsilon^{{1}/{2}}
\left(u_1\right)\right|^2\leq C\left(1+N^2\right)\left(u_2-u_1\right)^2 
\label{qq2}
\end{align}
for $\left|u_1\right|<N,\left|u_2\right|<N$
\end{lemma}
{\bf Proof.} The first estimate \eqref{qq1} we obtain immediately from
\eqref{qq}. The proof of the second estimate \eqref{qq2} can be carried out like 
the proof of the relation \eqref{rr}.

The properties of the process $\hat Z_\varepsilon \left(\cdot \right)$
established in the Lemmae 5 and 7 allows to cite Theorem 1.10.1 in \cite{IH81} and
to obtain the announced in the Theorem \ref{T2} properties of the pseudo-MLE
$\hat\vartheta _\varepsilon $.

\subsection{Discontinuous versus discontinuous} 

Let us remind that if the the observed model is discontinuous and the
statistician knows this but takes the wrong signals before and after the jump,
then nevereless it is possible to have the consistent estimation.  Consider
the following problem of parameter estimation in the situation of
misspecification.  The {\it theoretical model} is
\begin{align*}
{\rm d}X_t=\left[h\left(t\right)\1_{\left\{t<\vartheta
  \right\}}+g\left(t\right)\1_{\left\{t\geq \vartheta \right\}} \right]{\rm
  d}t+\varepsilon {\rm d}W_t,\quad X_0=0,\quad 0\leq t\leq T, 
\end{align*} 
where $\vartheta \in \Theta =\left(\alpha ,\beta \right)$, $0<\alpha <\beta
<T$. Suppose that $h\left(t\right)-g\left(t\right)>0$ for $t\in \left[\alpha
  ,\beta \right]$.  The observed stochastic process has  a different equation
\begin{align*}
{\rm
  d}X_t=\left[\left[h\left(t\right)+q\left(t\right)\right]\1_{\left\{t<\vartheta_0 
  \right\}}+\left[g\left(t\right)+r\left(t\right)\right]\1_{\left\{t\geq
    \vartheta_0 \right\}} \right]{\rm d}t+\varepsilon {\rm d}W_t,\quad 0\leq
t\leq T, 
\end{align*}
where $q\left(t\right)$ and $r\left(t\right)$ are some unknown functions. 

We study the conditions on $q\left(t\right)$ and $r\left(t\right)$ which allow
the consistent estimation of the parameter $\vartheta _0$.

The function $\Phi
\left(\vartheta \right)$ for $\vartheta <\vartheta _0$ is
\begin{align*}
\Phi \left(\vartheta \right)=\int_{0}^{\vartheta }q\left(t\right)^2{\rm
  d}t+\int_{\vartheta }^{\vartheta _0} \left[
  h\left(t\right)+q\left(t\right)-g\left(t\right)\right]^2 {\rm
  d}t+\int_{\vartheta _0}^{T} r\left(t\right)^2{\rm d}t.
\end{align*}
Hence
\begin{align*}
\frac{{\rm d}\Phi \left(\vartheta \right)}{{\rm d}\vartheta }
&=q\left(\vartheta \right)^2 -\left[h\left(\vartheta
  \right)-g\left(\vartheta\right)+q\left(\vartheta
  \right)\right]^2\\
&=-\left(h\left(\vartheta \right)-g\left(\vartheta
\right)\right) \left[h\left(\vartheta \right)-g\left(\vartheta
\right)+2q\left(\vartheta \right) \right].
\end{align*}
If the function 
\begin{equation}
\label{5-1}
q\left(\vartheta  \right)>\frac{g\left(\vartheta \right)-h\left(\vartheta
  \right)}{2} ,\qquad \vartheta \in\Theta ,
\end{equation}
then for $\vartheta <\vartheta _0$
\begin{align*}
\frac{{\rm d}\Phi \left(\vartheta \right)}{{\rm d}\vartheta }<0.
\end{align*}
For $\vartheta >\vartheta _0$ under condition 
\begin{equation}
\label{5-2}
r\left(\vartheta  \right)<\frac{h\left(\vartheta \right)-g\left(\vartheta
  \right)}{2} 
\end{equation}
we obtain the similar inequality 
\begin{align*}
\frac{{\rm d}\Phi \left(\vartheta \right)}{{\rm d}\vartheta }>0.
\end{align*}

Therefore 
\begin{align*}
\hat\vartheta =\arg\inf_{\vartheta \in \Theta }\Phi \left(\vartheta \right)=\vartheta _0
\end{align*}
and we obtain the following result.
\begin{proposition}
\label{P}
If the conditions \eqref{5-1} and \eqref{5-2} are  fulfilled then the pseudo-MLE
$\hat\vartheta _\varepsilon $ is consistent.
\end{proposition}

It can be shown that 
\begin{align*}
\frac{\hat\vartheta _\varepsilon-\vartheta _0}{\varepsilon ^2}\Longrightarrow \xi 
\end{align*}
For the details see the similar problem in Section 5.3, \cite{Kut94}. The
close problem of change-point detection formisspecified diffusion processes
are studied in \cite{CLGK00}.

\subsection{Discussion}

There are several other interesting problems of misspecification in
regularity, which can be studied by the proposed here approach. 

One of them is to study the asymptotic behavior of the bayesian estimator
$\tilde\vartheta _\varepsilon$  in 
the situation described by the equations \eqref{20}, \eqref{21}. The estimator
is 
\begin{align*}
\tilde\vartheta _\varepsilon =\frac{\int_{\alpha }^{\beta }\vartheta
  p\left(\vartheta \right)V\left(\vartheta ,X^T\right){\rm d}\vartheta
}{\int_{\alpha }^{\beta } p\left(\vartheta \right)V\left(\vartheta
  ,X^T\right){\rm d}\vartheta} ,
\end{align*}
where $p\left(\vartheta \right),\alpha <\vartheta <\beta $ is continuous
positive density of the  distribution of the random variable $\vartheta $.

It can be shown that $\tilde\vartheta _\varepsilon $ converges to the same
value $\hat\vartheta $. Then using the notations of the section 3.1 we can
write
\begin{align*}
\tilde\vartheta _\varepsilon& =\frac{\int_{\alpha }^{\beta }\vartheta
  p\left(\vartheta \right) \frac{V\left(\vartheta ,X^T\right)}{
    V\left(\hat\vartheta ,X^T\right)}{\rm d}\vartheta 
}{\int_{\alpha }^{\beta } p\left(\vartheta \right)\frac{V\left(\vartheta ,X^T\right)}{
    V\left(\hat\vartheta ,X^T\right)}{\rm d}\vartheta} =\hat\vartheta
+\varepsilon ^{2/3} \frac{\int_{\UU_\varepsilon }^{}u\,p\left(\vartheta
  _u\right)Z_\varepsilon \left(u\right){\rm d}u}{\int_{\UU_\varepsilon }^{}p\left(\vartheta
  _u\right)Z_\varepsilon \left(u\right){\rm d}u} ,
\end{align*}
where we changed the variables $\vartheta =\vartheta _u=\hat\vartheta
+\varepsilon ^{2/3}u $. Hence 
\begin{align*}
\frac{\tilde\vartheta _\varepsilon-\hat\vartheta  }{\varepsilon ^{2/3}}
\approx \frac{\int_{\UU_\varepsilon }^{}u\,Z_\varepsilon \left(u\right){\rm
    d}u}{\int_{\UU_\varepsilon }^{}\,Z_\varepsilon \left(u\right){\rm
    d}u}=\frac{\int_{\UU_\varepsilon }^{}u\, \left(\hat Z_\varepsilon
  \left(u\right)\right)^{2\varepsilon ^{-2/3}}  {\rm
    d}u}{\int_{\UU_\varepsilon }^{}\,\left(\hat Z_\varepsilon
  \left(u\right)\right)^{2\varepsilon ^{-2/3}}{\rm d}u}
\end{align*}
and the problem reduces to the study of the asymptotics of these two
integrals in the situation, when
\begin{align*}
\hat Z_\varepsilon \left(u\right)=\exp\left\{\delta (\hat\vartheta
)W\left(u\right) -\frac{\gamma (\hat\vartheta
  )}{2}u^2\right\} \left(1+o\left(1\right)\right).
\end{align*}
We can suppose that the detailed study will provide us the asymptotics 
\begin{align*}
\frac{\tilde\vartheta _\varepsilon-\hat\vartheta  }{\varepsilon ^{2/3}}
\approx  \hat u
\end{align*} 
where $\hat u $ is as before the point of the maximum of the process $\delta (\hat\vartheta
)W\left(u\right) -\frac{\gamma (\hat\vartheta   )}{2}u^2$. This means that as
usual in regular estimation problems the asymptotic behavior of the  BE
is similar to that of the MLE. 

Another problem we obtain  if we suppose that the obsevred process has a
signal $M\left(\vartheta ,\cdot \right)$ with a singularity of the
{\it cusp-type} (theoretical model) but the observed process in reality has a
{\it smooth} signal $S\left(\vartheta ,\cdot \right)$, i.e.; 
{\it cusp versus smooth}.  Say, 
\begin{align}
\label{cusp}
{\rm d}X_t=a\left|t-\vartheta \right|^\kappa {\rm d}t+\varepsilon {\rm
  d}W_t,\quad X_0=0,\quad 0\leq t\leq T ,
\end{align}
where $\kappa \in (0,\frac{1}{2})$. 

Therefore the observed process is \eqref{1-6} but the statistician calculate
the LR $V\left(\vartheta ,X^T\right)$ and the pseudo-MLE $\hat\vartheta
_\varepsilon $ following \eqref{12} and \eqref{13} respectively. It is clear
that $\hat\vartheta
_\varepsilon  $ converges to the value
\begin{align*}
\hat\vartheta =\arg\inf_{\vartheta \in \Theta
}\int_{0}^{T}\left[a\left|t-\vartheta \right|^\kappa-S\left(\vartheta
  _0,t\right) \right]^2{\rm d}t 
\end{align*}
which minimizes the Kullback-Leibner distance and we are interested by the
limit distribution of  $\varepsilon 
^{-\frac{2}{3-2\kappa }}\left(\hat\vartheta 
_\varepsilon -\hat\vartheta\right)$. For the details see the forthcomming work
\cite{CDK15}.

 The properties of the MLE and bayesian
estimatrs for the ergodic diffusion processes and inhomogeneous Poisson
processes 
with cusp-type singularities  are studied for example, in \cite{DK03}, \cite{D03}. 
For the general theory of the parameter estimation for different singular
estimation problems see \cite{IH81}.
 
\bigskip

{\bf Acknowledgment.} This work was done under partial financial support of
the grant of  RSF number 15-11-10022.

\end{document}